\documentclass[12pt]{article}

\usepackage[cp1251]{inputenc}
\usepackage[english,russian,ukrainian]{babel}

\usepackage{latexsym,upref,amsfonts,amsthm,cite,amssymb,amsmath,epsfig}
\usepackage{graphicx,graphics,hhline}
\usepackage{euscript}

\begin{document}

\theoremstyle{plain}
\newtheorem{theorem}{Теорема}
\newtheorem{lemma}{Лема}
\newtheorem{corollary}{Наслідок}
\newtheorem{property}{Властивість}
\theoremstyle{definition}
\newtheorem{definition}{Означення}
\theoremstyle{remark}
\newtheorem{remark}{Зауваження}

\selectlanguage{ukrainian}

\noindent {\small УДК 517.5} \vskip 3mm

\noindent \textbf{A.S. Serdyuk, I.V. Sokolenko} \\ {\small (Institute of Mathematics of National Academy of Sciences of Ukraine, Kiev)}\vskip 5mm

\noindent \textbf{А.\,С. Сердюк, І.\,В. Соколенко} \\  {\small (Iнститут математики Національної академії наук України, Київ)} \vskip 10mm

 \noindent  \textbf{\Large Approximation  of classes of convolutions of \\ periodic functions by linear methods constructed \\ on basis of Fourier-Lagrange coefficients}
 \vskip 10mm

\noindent  \textbf{\Large Апроксимація  класів згорток періодичних функцій \\   лінійними методами, побудованими \\ на основі   коефіцієнтів Фур'є--Лагранжа}
 \vskip 10mm

{\small \it \noindent We calculate the least upper bounds of pointwise and uniform approximations for classes  of $2\pi$-periodic functions expressible as convolutions of an arbitrary square summable  kernel with functions, which  belong to  the unit ball of the space $L_2$,   by linear polynomial methods  constructed on the basis of their Fourier-Lagrange coefficients.
 \hfill}
 \vskip 7mm

{\small \it \noindent Знайдено  точні верхні межі  поточкових та рівномірних наближень  класів $2\pi$-періодичних функцій, що задаються у вигляді згортки довільного сумовного з квадратом твірного ядра з функціями, що належать одиничній кулі  простору $L_2$, за допомогою лінійних поліноміальних методів, побудованих на основі їх коефіцієнтів Фур'є--Лагранжа. \hfill}

\newpage
\normalsize\vskip 3mm

Нехай $C$ i
$L_p, \ 1\le p\le\infty,$ --- простори
$2\pi$-періодичних  функцій зі стандартними нормами $\|\cdot\|_C$ та $\|\cdot\|_p,$ відповідно.

Позначимо через  $C^\psi_{\bar\beta,p},\ 1\le p\le\infty,$  множину всіх $2\pi$-періодичних функцій $f$, які зображуються за допомогою  згортки
\begin{equation}\label{7'}
    f(x)=\frac{a_0}{2}+\frac{1}{\pi}\int\limits_{-\pi}^{\pi}
 \varphi(x-t) \Psi_{\bar\beta}(t)dt, \ \ \ a_0\in\mathbb R, \ \ \ \varphi\in B_p^0,
\end{equation}

$$
  B_p^0=\{\varphi\in L_p:\ \|\varphi\|_p\le1,\ \varphi\perp1\},
$$
із фіксованим твірним ядром $\Psi_{\bar\beta}\in L_{p'},\  1/p+ 1/{p'}=1,\ $ ряд Фур'є якого має вигляд
\begin{equation*}\label{3}
 \Psi_{\bar\beta}(t)\sim\sum\limits_{k=1}^\infty \psi(k)\cos\left(kt-\frac{\beta_k\pi}2\right),
\end{equation*}
де $\psi=\psi(k)$ i $\bar\beta=\beta_k$, $k=1,2,\ldots,$ --- довільні послідовності
дійсних чисел. Якщо $\psi(k)\neq0,\ k\in\mathbb{N},$ то
функцію $\varphi$ у зображенні (\ref{7'}) називають $(\psi,\bar{\beta})$-похідною функції $f$ і позначають  $f^{ \psi}_{\bar{\beta}}$. Поняття $(\psi,\bar{\beta})$-похідної введено О.~І.~Степанцем (див., наприклад, \cite{Stepanets2002_1}).  Оскільки \mbox{$\varphi\in L_p,$} a $\ \Psi_{\bar\beta}\in L_{p'},$ то (див. \cite[c.~144]{Stepanets2002_1}) згортка (\ref{7'}) є неперервною функцією, тобто $C^\psi_{\bar\beta,p}\subset C.$

При $\beta_k\equiv\beta,\ \beta\in \mathbb R,$ класи $C^\psi_{\bar\beta,p}$ позначатимемо через $C^\psi_{\beta,p}$. При $\psi(k)=k^{-r}, r>0, $ класи $C^\psi_{\bar\beta,p}$ будемо позначати через $W^r_{\bar\beta,p}$. Якщо \mbox{$\psi(k)=k^{-r}$} i  $\beta_k\equiv r,\ r\in\mathbb N,$ то   $W^r_{\bar\beta,p}$ є відомими класами $r-$диференційовних функцій $W^r_p$. Якщо  $\psi(k)=q^{k},$  \mbox{$0< q<1,$} то  класи  $C^\psi_{\bar\beta,p}$ будемо позначати через $C^q_{\bar\beta,p}.$  При $\beta_k\equiv\beta,\ \beta\in \mathbb R,$ $C^q_{\bar\beta,p}$ є відомими класами інтегралів Пуассона  $C^q_{\beta,p}.$
Зрозуміло, що при $p=2$ умова включення $\Psi_{\bar\beta}\in L_2$ еквівалентна виконанню умови
\begin{equation}\label{6}
    \sum_{k=1}^\infty \psi^2(k)<\infty.
\end{equation}

Нехай $\ f(x)$ --- довільна $\ 2\pi$-періодична неперервна
функція. Через $\ \widetilde  S_n(f; x)\ $будемо позначати
тригонометричний поліном порядку $\ n,\ $що інтерполює $\ f(x)\ $у точках $\ x_k^{(n)}=2k\pi/(2n+1),\ $ \mbox{$k=0,1,\ldots,
2n,$} тобто такий, що $$\widetilde  S_n(f; x_k^{(n)}) = f(x_k^{(n)}), \
\ k=0, 1, \ldots, 2n.$$

 Інтерполяційний тригонометричний поліном $\widetilde{S}_n(f;x)$
можна записати (див., наприклад, \cite[c. 128-129]{Stepanets2002_2}) в явному вигляді наступним чином:
\begin{equation}\label{10}
 \widetilde{S}_n(f;x) =\frac{a_0^{(n)}}{2}+\sum\limits_{k=1}^n (a_k^{(n)}
 \cos kx+b_k^{(n)}\sin kx),
\end{equation}
 де
 $$
 a_k^{(n)}=\frac{2}{2n+1}\sum\limits_{i=0}^{2n} f(x_i^{(n)}) \cos
 kx_i^{(n)},\quad k=0,1,2,\ldots ,
 $$
 $$
 b_k^{(n)}=\frac{2}{2n+1}\sum\limits_{i=0}^{2n} f(x_i^{(n)}) \sin
 kx_i^{(n)},\quad k=1,2,\ldots.
  $$

 Зв'язок між коефіцієнтами Фур'є $a_k$ i $b_k$ функції $f(x)$ \mbox{$(f\in
 C)$} i коефіцієнтами $a_k^{(n)}$ i $b_k^{(n)}$ інтерполяційного
 полінома $\widetilde{S}_n(f;x)$ виражається за допомогою
 рівностей \cite[c. 130]{Stepanets2002_2}
 \begin{equation}\label{11}
   a_k^{(n)}=a_k+\sum\limits_{m=1}^\infty
  (a_{m(2n+1)+k}+a_{m(2n+1)-k}),
  \quad k=0,1,2,\ldots ,
 \end{equation}
\begin{equation}\label{12}
  b_k^{(n)}=b_k+\sum\limits_{m=1}^\infty
  (b_{m(2n+1)+k}-b_{m(2n+1)-k}),
  \quad k=1,2,\ldots
\end{equation}

Розглянемо лінійні поліноміальні методи наближення  функцій $f$ з класів $C^\psi_{\bar\beta,2}$, які  побудовані на основі їх коефіцієнтів Фур'є--Лагранжа $a_k^{(n)}$ i $b_k^{(n)}$. Нехай $\Lambda=\|\lambda_k^{(n)}\|\ $ i $\ {\rm M}=\|\mu_k^{(n)}\|,$  $n=0,1,\ldots,  $ \mbox{$\ k=0,1,\ldots,$~---} нескінченні трикутні матриці дійсних чисел такі, що
\begin{equation}\label{a1}
   \begin{array}{lll}
    \lambda_0^{(n)}=1,& \mu_0^{(n)}=0, &  n=0,1,2,\ldots,\\ \\
    \lambda_k^{(n)}=0, &  \mu_k^{(n)}=0, & k=n+1, n+2, \ldots, \\ \\
    \lim\limits_{n\rightarrow \infty}\lambda_k^{(n)}=1,& \lim\limits_{n\rightarrow \infty}\mu_k^{(n)}=0, & k=1,2,\ldots
  \end{array}
\end{equation}
Позначимо через  $\widetilde{U}_n=\widetilde{U}_n(\Lambda;{\rm M})$ лінійний оператор, який кожній функції $f\in C$ 
ставить у відповідність тригонометричний поліном вигляду
\begin{equation}\label{a2}
     \widetilde{U}_n(f;x)=\widetilde{U}_n(f;\Lambda;{\rm M};x)=
\frac{a_0^{(n)}}{2}+\sum\limits_{k=1}^n\left(\lambda_k^{(n)}(a_k^{(n)}\cos kx+b_k^{(n)}\sin kx)+\right.
$$
$$
\left.+\mu_k^{(n)}(-b_k^{(n)}\cos kx+a_k^{(n)}\sin kx)\right),
\end{equation}
де $a_k^{(n)}$\ i\ $b_k^{(n)}$ --- коефіцієнти, що визначаються рівностями  $(\ref{11})$ i $(\ref{12})$.

У даній роботі розглядається задача про знаходження точних значень величин
  \begin{equation}\label{9''}
    \widetilde{\cal E}_n(C^\psi_{\bar\beta,2};\Lambda;{\rm M};x)=\sup\limits_{f\in C^\psi_{\bar\beta,2}}|f(x)-\widetilde  U_n(f;x)|.
 \end{equation}
 Має місце наступне твердження.

\vspace*{2mm}\begin{theorem}\label{1t}
 Нехай послідовність дійсних чисел $\psi(k)$ задовольняє умову $(\ref{6})$,  а  \mbox{$\Lambda=\|\lambda_k^{(n)}\|\ $} i $\ {\rm M}=\|\mu_k^{(n)}\|$ ---  умови $(\ref{a1})$. Тоді для довільної послідовності $\bar\beta=\beta_k, \beta_k\in\mathbb R,$ і довільного $n\in \mathbb N$ у кожній точці $x\in\mathbb R$ виконується рівність
$$
  \widetilde{\cal E}_n(C^\psi_{\bar\beta,2};\Lambda;{\rm M};x)=
 \frac1{\sqrt\pi}\Bigg(\sum_{k=1}^n\left((1-\lambda_k^{(n)})^2+(\mu_k^{(n)})^2\right)\psi^2(k)+
$$
$$
+ \sum\limits_{m=1}^\infty\sum_{k=m(2n+1)-n}^{m(2n+1)+n}\left(\left(\cos m(2n+1)x -\lambda_{|k-m(2n+1)|}^{(n)}\right)^2+\right.
$$
$$
+
\left.\left(\sin m(2n+1)x +\mu_{|k-m(2n+1)|}^{(n)}\right)^2
\right)\psi^2(k)\Bigg)^{1/2}.
$$
 \end{theorem}

\noindent\textit{\textbf{Доведення.}} Для довільної функції $  f$ з класу $C^\psi_{\bar\beta,2}$  згідно з $(\ref{7'})$ мають місце формули
\begin{equation}\label{7_1}
   a_\nu\cos kx+b_\nu\sin kx=\frac1\pi\int\limits_{-\pi}^\pi \psi(\nu)\cos\left(\nu(x-t)-(\nu-k)x-\frac{\beta_\nu\pi}2\right)\varphi(t)dt,
\end{equation}
\begin{equation}\label{7_2}
 -b_\nu\cos kx+a_\nu\sin
 kx=\frac1\pi\int\limits_{-\pi}^\pi \psi(\nu)\sin\left(\nu(x-t)-(\nu-k)x-\frac{\beta_\nu\pi}2\right)\varphi(t)dt.
\end{equation}

 Об'єднуючи (\ref{7'}), (\ref{11})--(\ref{a2}), (\ref{7_1}) та (\ref{7_2}), для $f\in C^\psi_{\bar\beta,2}$ одержуємо
\begin{equation}\label{l}
f(x)-\widetilde  U_n(f;x)=
$$
$$
=\frac 1{\pi }\int\limits_{-\pi }^{\pi }\varphi(x-t)
\Bigg(\sum\limits_{k=1}^n\psi(k)\left((1-\lambda_k^{(n)})\cos\left(k t-\frac{\beta_k\pi}2\right)-\mu_k^{(n)}\sin\left(k t-\frac{\beta_k\pi}2\right)\right)+
$$
$$
  +
\sum _{m=1}^{\infty }  \sum _{\nu =m(2n+1)-n}^{m(2n+1)+n}
\psi (\nu )
\left(\cos \left(\nu t-\frac {\beta_\nu \pi}2\right)-\right.
$$
$$
-
\lambda_{|\nu-m(2n+1)|}^{(n)}\cos\left(\nu t-m(2n+1)x-\frac{\beta_\nu\pi}2\right)
-
$$
$$
\left.
-\mu_{|\nu-m(2n+1)|}^{(n)}\sin\left(\nu t-m(2n+1)x-\frac{\beta_\nu\pi}2\right) \right)\Bigg)dt=
$$
$$
=
\frac 1{\pi }\int\limits_{-\pi }^{\pi }\varphi(x-t)
\Bigg(\sum\limits_{k=1}^n\psi(k)\left((1-\lambda_k^{(n)})\cos\left(k t-\frac{\beta_k\pi}2\right)-\right.
\left.\mu_k^{(n)}\sin\left(k t-\frac{\beta_k\pi}2\right)\right)+
$$
$$
+
\sum _{m=1}^{\infty }  \sum _{\nu =m(2n+1)-n}^{m(2n+1)+n}
\psi (\nu )
\left(\left(\cos m(2n+1)x -\lambda_{|\nu-m(2n+1)|}^{(n)}\right) \cos \left(\nu t-m(2n+1)x-\frac{\beta_\nu\pi}2\right)-\right.
$$
$$
\left.
-\left(\sin m(2n+1)x +\mu_{|\nu-m(2n+1)|}^{(n)}\right)\sin\left(\nu t-m(2n+1)x-\frac{\beta_\nu\pi}2\right) \right)\Bigg)dt.
\end{equation}

Розглядаючи точну верхню межу по функціях $f$ з класу $C^\psi_{\bar\beta,2}$ у лівій і правій частині  рівності (\ref{l})  та враховуючи інваріантність множин $C^\psi_{\bar\beta,2}$ відносно зсуву аргументу,  в кожній точці $x\in\mathbb R$ отримуємо
 \begin{equation}\label{14}
    \widetilde{\cal E}_n(C^\psi_{\bar\beta,2};\Lambda;{\rm M};x)=
    $$
    $$
    =
\frac 1{\pi} \sup\limits_{\varphi\in B_2^0}\int\limits_{-\pi }^{\pi }\varphi(t)
\Bigg(\sum\limits_{k=1}^n\psi(k)\left((1-\lambda_k^{(n)})\cos\left(k t-\frac{\beta_k\pi}2\right)-\mu_k^{(n)}\sin\left(k t-\frac{\beta_k\pi}2\right)\right)+
$$
$$
+\sum _{m=1}^{\infty }  \sum _{\nu =m(2n+1)-n}^{m(2n+1)+n}
\psi (\nu )
\left(\left(\cos m(2n+1)x -\lambda_{|\nu-m(2n+1)|}^{(n)}\right) \cos \left(\nu t-m(2n+1)x-\frac{\beta_\nu\pi}2\right)-\right.
$$
$$
\left.
-\left(\sin m(2n+1)x +\mu_{|\nu-m(2n+1)|}^{(n)}\right)\sin\left(\nu t-m(2n+1)x-\frac{\beta_\nu\pi}2\right) \right)\Bigg)dt.
\end{equation}

Застосовуючи до правої частини (\ref{14})  співвідношення двоїстості \cite[c.~27]{Kornejchuk1987}
$$
 \sup_{\varphi\in B_2^0} \int\limits_{-\pi}^\pi \varphi(t)u(t)dt=\inf_{\alpha\in \mathbb R}\|u-\alpha \|_2,\ \ \ u\in L_2,
 $$
  та рівність Парсеваля, одержуємо
   $$
    \widetilde{\cal E}_n(C^\psi_{\bar\beta,2};\Lambda;{\rm M};x)=
    \frac1\pi\inf_{\alpha\in \mathbb R}\Bigg\|\sum\limits_{k=1}^n\psi(k)\left((1-\lambda_k^{(n)})\cos\left(k t-\frac{\beta_k\pi}2\right)-\right.
\left.\mu_k^{(n)}\sin\left(k t-\frac{\beta_k\pi}2\right)\right)
+
$$
$$
+\sum _{m=1}^{\infty }  \sum _{\nu =m(2n+1)-n}^{m(2n+1)+n}
\psi (\nu )
\left(\left(\cos m(2n+1)x -\lambda_{|\nu-m(2n+1)|}^{(n)}\right) \cos \left(\nu t-m(2n+1)x-\frac{\beta_\nu\pi}2\right)-\right.
$$
$$
\left.
-\left(\sin m(2n+1)x +\mu_{|\nu-m(2n+1)|}^{(n)}\right)\sin\left(\nu t-m(2n+1)x-\frac{\beta_\nu\pi}2\right) \right)-\alpha\Bigg\|_2=
   $$
   $$
   =\frac1{\sqrt\pi}\inf_{\alpha\in \mathbb R}\Bigg(\sum_{k=1}^n\left((1-\lambda_k^{(n)})^2+(\mu_k^{(n)})^2\right)\psi^2(k)+
$$
$$
+ \sum\limits_{m=1}^\infty\sum_{k=m(2n+1)-n}^{m(2n+1)+n}\left(\left(\cos m(2n+1)x -\lambda_{|k-m(2n+1)|}^{(n)}\right)^2+\right.
$$
$$
+
\left.\left(\sin m(2n+1)x +\mu_{|k-m(2n+1)|}^{(n)}\right)^2
\right)\psi^2(k)+\alpha^2\Bigg)^{1/2}=
$$
$$
=
 \frac1{\sqrt\pi}\Bigg(\sum_{k=1}^n\left((1-\lambda_k^{(n)})^2+(\mu_k^{(n)})^2\right)\psi^2(k)+
$$
$$
+ \sum\limits_{m=1}^\infty\sum_{k=m(2n+1)-n}^{m(2n+1)+n}\left(\left(\cos m(2n+1)x -\lambda_{|k-m(2n+1)|}^{(n)}\right)^2+\right.
$$
$$
+
\left.\left(\sin m(2n+1)x +\mu_{|k-m(2n+1)|}^{(n)}\right)^2
\right)\psi^2(k)\Bigg)^{1/2}.
$$
Теорему 1 доведено.

Якщо
\begin{equation*}\label{n1}
    \begin{array}{ccc}
      \lambda_k^{(n)}=\left\{
                   \begin{array}{ll}
                     1, & 0\le k\le n, \\
                     0, & k>n,
                   \end{array}
                 \right. & i & \mu_k^{(n)}\equiv0,
    \end{array}
\end{equation*}
то згідно з (\ref{10}) і (\ref{a2}) тригонометричний поліном $\widetilde U_n(f;\Lambda;{\rm M};x)$ є інтерполяційним тригонометричним поліномом $\widetilde{S}_n(f;x)$ функції $f$ порядку $n$, тобто $\widetilde U_n(f;\Lambda;{\rm M};x)=\widetilde{S}_n(f;x)$. В цьому випадку, як випливає з теореми 1, при всіх  $x\in\mathbb{R}$ для величин
\begin{equation}\label{15'}
{\cal E}(C^\psi_{\bar\beta,2}; \widetilde S_n;x) =\sup\limits_{f\in C^\psi_{\bar\beta,2}}|f(x)-\widetilde  S_n(f;x)|,\ \ \ n\in\mathbb{N},
\end{equation}
виконується рівність
 \begin{equation}\label{15}
 {\cal E}(C^\psi_{\bar\beta,2}; \widetilde S_n;x) =
  \frac2{\sqrt\pi}\left(\sum\limits_{m=1}^\infty\sin^2\frac{(2n{+}1)mx}2\sum_{k=m(2n+1)-n}^{m(2n+1)+n}\psi^2(k)\right)^{1/2}.
 \end{equation}
Рівність (\ref{15})  встановлено в \cite{Serdyuk_Sokolenko2016}.

Наряду з  ${\cal E}(C^\psi_{\bar\beta,2}; \widetilde S_n;x)$ розглянемо величини
\begin{equation}\label{15''}
{\cal E}(C^\psi_{\bar\beta,2}; \widetilde S_n)_C =
\sup\limits_{f\in C^\psi_{\bar\beta,2}}\|f(\cdot)-\widetilde  S_n(f;\cdot)\|_C, \ \ \ n\in\mathbb{N}.
\end{equation}
Має місце наступне твердження.

\vspace*{2mm}\begin{theorem}\label{2t}
 Нехай $n\in \mathbb N$, а послідовність дійсних чисел $\psi(k)$ задовольняє умову $(\ref{6})$ і така, що послідовність $\alpha_m=m\sum\limits_{k=m(2n+1)-n}^{m(2n+1)+n}\psi^2(k)$ є опуклою донизу. Тоді для довільної  $\bar\beta=\beta_k, \beta_k\in\mathbb R,$ виконуються рівності
 $$
 {\cal E}(C^\psi_{\bar\beta,2}; \widetilde S_n)_C={\cal E}(C^\psi_{\bar\beta,2}; \widetilde S_n;\frac\pi{2n{+}1}) =
\frac2{\sqrt\pi}\left(\sum\limits_{l=1}^\infty \sum_{k=(2l-1)(2n+1)-n}^{(2l-1)(2n+1)+n}\psi^2(k)\right)^{1/2}.
 $$
 \end{theorem}

\noindent\textit{\textbf{Доведення.}} Підносячи до квадрату обидві частини рівності (\ref{15}), отримуємо
$$
 {\cal E}^2(C^\psi_{\bar\beta,2}; \widetilde S_n;x) =
  \frac2{ \pi}\sum\limits_{m=1}^\infty(1-\cos(2n+1)mx)\sum_{k=m(2n+1)-n}^{m(2n+1)+n}\psi^2(k).
$$
При кожному $n\in\mathbb N$ величина ${\cal E}^2(C^\psi_{\bar\beta,2}; \widetilde S_n;x)$  є парною періодичною функцією від змінної $x$ з періодом $\displaystyle T= \frac{2\pi}{2n+1},$ тому  досить її розглядати на
$\displaystyle [0,\frac\pi{2n{+}1}].$ За виконання умов теореми 2, як випливає з  \cite[c.~297]{Zygmund1965}, похідна
$$
\frac{d}{dx}{\cal E}^2(C^\psi_{\bar\beta,2}; \widetilde S_n;x)=\frac{2(2n+1)}{\pi}\sum\limits_{m=1}^\infty m\sum_{k=m(2n+1)-n}^{m(2n+1)+n}\psi^2(k) \sin(2n+1)mx
$$
 є додатною на $(\displaystyle 0,\frac\pi{2n{+}1})$  функцією, а, отже,  максимум величини ${\cal E}^2(C^\psi_{\bar\beta,2}; \widetilde S_n;x)$
 досягається в точці $\displaystyle x=\frac\pi{2n{+}1}$ i приймає значення
$$
  {\cal E}^2(C^\psi_{\bar\beta,2}; \widetilde S_n)_C=\max_{x\in \mathbb R}{\cal E}^2(C^\psi_{\bar\beta,2}; \widetilde S_n;x)={\cal E}^2(C^\psi_{\bar\beta,2}; \widetilde S_n; \frac\pi{2n{+}1}) =
  $$
  $$
  =\frac4{ \pi}\sum\limits_{l=1}^\infty  \sin^2\frac{(2l{-}1)\pi}2\sum_{k=(2l-1)(2n+1)-n}^{(2l-1)(2n+1)+n}\psi^2(k)
  +\sin^2l\pi \sum_{k=2l(2n+1)-n}^{2l(2n+1)+n}\psi^2(k)=
   $$
  $$
  =
   \frac4{ \pi}\sum\limits_{l=1}^\infty \sum_{k=(2l-1)(2n+1)-n}^{(2l-1)(2n+1)+n}\psi^2(k).
$$
Теорему 2 доведено.

\noindent\textit{\textbf{Зауваження.}} Із означень (\ref{15'}) i (\ref{15''}) випливають очевидні співвідношення
$$
{\cal E}(C^\psi_{\bar\beta,2}; \widetilde S_n;x)\le \max_{x\in \mathbb R}{\cal E}(C^\psi_{\bar\beta,2}; \widetilde S_n;x)={\cal E}(C^\psi_{\bar\beta,2}; \widetilde S_n)_C. $$
Як випливає з теореми  2, опуклість донизу послідовності \mbox{$\alpha_m=m\sum\limits_{k=m(2n+1)-n}^{m(2n+1)+n}\psi^2(k)$} є достатньою умовою того, що величини ${\cal E}(C^\psi_{\bar\beta,2}; \widetilde S_n;x)$ (\ref{15'}) досягають значень ${\cal E}(C^\psi_{\bar\beta,2}; \widetilde S_n)_C$ рівно посередині між вузлами інтерполяції, тобто при $\displaystyle x=\frac\pi{2n{+}1} +\frac{2j\pi}{2n{+}1}, j\in\mathbb Z.$
Неважко переконатись, що при достатньо великих $n$ умову опуклості послідовності $\alpha_m$ задовольняють  \mbox{$\psi(k)=q^k,$}  $q\in(0,1).$
Утім,  наступне твердження показує, що для зазначених $\psi(k)$ рівності
 \begin{equation}\label{1z1}
{\cal E}(C^\psi_{\bar\beta,2}; \widetilde S_n; \frac\pi{2n{+}1})={\cal E}(C^\psi_{\bar\beta,2}; \widetilde S_n)_C,
\end{equation}
мають місце при всіх $n\in\mathbb N$.

\vspace*{2mm} \begin{theorem}\label{3t}
 Нехай $q\in(0,1),\ $ $\bar\beta=\beta_k$~---  довільна послідовність дійсних чисел  і  $n\in \mathbb N$. Тоді
$$
 {\cal E}(C^q_{\bar\beta,2}; \widetilde S_n)_C={\cal E}(C^q_{\bar\beta,2}; \widetilde S_n;\frac\pi{2n{+}1}) =\frac{2q^{n+1}}{\sqrt{\pi(1-q^2)(1+q^{2(2n+1)})}}.
$$
 \end{theorem}

\noindent\textit{\textbf{Доведення.}} Згідно з наслідком 1 роботи \cite{Serdyuk_Sokolenko2016}
\begin{equation}\label{3t1}
 {\cal E}(C^q_{\bar\beta,2}; \tilde S_n;x) =\left|\sin\frac{(2n+1)x}2\right|\frac{2q^{n+1}}{\sqrt{\pi(1-q^2)}}\times
 $$
 $$
 \times\left(\frac{1+q^{2(2n+1)}}{1-2q^{2(2n+1)}\cos(2n+1)x+q^{4(2n+1)}}\right)^{1/2},\ \ \ x\in\mathbb R.
 \end{equation}
Рівність (\ref{3t1}) можна записати в еквівалентному вигляді:
 \begin{equation}\label{3t2}
 {\cal E}^2(C^q_{\bar\beta,2}; \tilde S_n;x)=\frac{2q^{2(n+1)}(1+q^{2(2n+1)})(1-\cos(2n+1)x)}{\pi(1-q^2)(1-2q^{2(2n+1)}\cos(2n+1)x+q^{4(2n+1)})}.
 \end{equation}
Розглянемо функцію
$$
 y(\tau)=\frac{1-\cos\tau}{1-2\rho\cos\tau+\rho^2},\ \ \ \rho\in(0,1).
$$
Оскільки
$$
y'(\tau)=\frac{(1-\rho)^2\sin\tau}{(1-2\rho\cos\tau+\rho^2)^2},
$$
то  максимум функції $y(\tau)$ досягається в точці $\tau=\pi$ і приймає значення
$$
y_{\max}=y(\pi)=\frac2{(1+\rho)^2}.
$$

Отже, з урахуванням (\ref{3t2}),
$$
 {\cal E}^2(C^q_{\bar\beta,2}; \tilde S_n)_C=
 \max_{x\in\mathbb R}{\cal E}^2(C^q_{\bar\beta,2}; \tilde S_n;x)=
  \frac{4q^{2(n+1)}(1+q^{2(2n+1)})}{\pi(1-q^2)(1+q^{2(2n+1)})^2}=\frac{4q^{2(n+1)}}{\pi(1-q^2)(1+q^{2(2n+1)})}.
 $$
Теорему 3 доведено.

Наступне твердження показує, що (\ref{1z1}) має місце для \mbox{$\psi(k)=k^{-r},$} $r> 1/2.$

\vspace*{2mm}\begin{theorem}\label{4t}
 Нехай $\psi(k)=k^{-r},\ r> 1/2,\ $ $\bar\beta=\beta_k$~---  довільна послідовність дійсних чисел  і  $n\in \mathbb N$. Тоді
 $$
 {\cal E}(W^r_{\bar\beta,2}; \widetilde S_n)_C={\cal E}(C^\psi_{\bar\beta,2}; \widetilde S_n;\frac\pi{2n{+}1}) =
 \frac2{\sqrt{\pi\Gamma(2r)}(2n+1)^{r}}
 \left(\int\limits_0^1\frac{\rho^{-n/(2n+1)}\ln^{2r-1}\displaystyle \rho^{-1} 
 }{(1-\rho^{1/(2n+1)})(1+\rho)}d\rho \right)^{1/2},
$$
 де $\Gamma(x)$ --- гамма-функція.
\end{theorem}

\noindent\textit{\textbf{Доведення теореми 4.}} Згідно з наслідком 2  роботи \cite{Serdyuk_Sokolenko2016}
 \begin{equation}\label{4t1}
 {\cal E}(W^r_{\bar\beta,2}; \tilde S_n;x) =\left|\sin\frac{(2n+1)x}2\right|\frac2{\sqrt{\pi\Gamma(2r)}(2n+1)^{r}}\times
 $$
 $$
 \times\left(\int\limits_0^1\frac{\rho^{-n/(2n+1)}(1+\rho)\ln^{2r-1}\displaystyle \rho^{-1} 
 }{(1-\rho^{1/(2n+1)})(1-2\rho\cos(2n+1)x+\rho^2)}d\rho \right)^{1/2}, \ x\in\mathbb R.
\end{equation}
Запишемо рівність (\ref{4t1})   у наступному вигляді:
 \begin{equation}\label{4t2}
 {\cal E}^2(W^r_{\bar\beta,2}; \tilde S_n;x)=\frac{2(1-\cos(2n+1)x)}{{\pi\Gamma(2r)(2n+1)^{2r}}}\times
 $$
 $$
 \times\int\limits_0^1\frac{\rho^{-n/(2n+1)}(1+\rho)\ln^{2r-1}\displaystyle \rho^{-1} 
 }{(1-\rho^{1/(2n+1)})(1-2\rho\cos(2n+1)x+\rho^2)}d\rho.
 \end{equation}
З рівності (\ref{4t2}) отримуємо
\begin{equation}\label{4t3}
 \frac{d}{dx}{\cal E}^2(W^r_{\bar\beta,2}; \widetilde S_n;x)=
$$  $$
 =\frac{2(2n+1)\sin(2n+1)x}{{\pi\Gamma(2r)(2n+1)^{2r}}}
 \int\limits_0^1\frac{\rho^{-n/(2n+1)}(1+\rho)(1-\rho)^2\ln^{2r-1}\displaystyle \rho^{-1} 
 }{(1-\rho^{1/(2n+1)})(1-2\rho\cos(2n+1)x+\rho^2)^2}d\rho.
 \end{equation}
Оскільки інтеграл в правій частині (\ref{4t3}) додатній, то  максимум ${\cal E}^2(W^r_{\bar\beta,2}; \tilde S_n;x)$ досягається в точці $\displaystyle x=\frac\pi{2n{+}1}$ і, з урахуванням (\ref{4t2}), приймає значення
$$
 {\cal E}^2(W^r_{\bar\beta,2}; \widetilde S_n)_C=
 \max_{x\in\mathbb R}{\cal E}^2(W^r_{\bar\beta,2}; \tilde S_n;x)=
 \frac4{\pi\Gamma(2r)(2n+1)^{2r}}
  \int\limits_0^1\frac{\rho^{-n/(2n+1)}\ln^{2r-1}\displaystyle \rho^{-1} 
 }{(1-\rho^{1/(2n+1)})(1+\rho)}d\rho.
$$
Теорему 4 доведено.


\bigskip

\small{

}

\vskip 10mm
E-mail: serdyuk@imath.kiev.ua,  sokol@imath.kiev.ua

 \end{document}